\documentclass[a4paper,12pt,draft]{article}
\usepackage{amsmath,amssymb}
\title{Natural volume forms on pseudo-Finslerian manifolds with $m$-th root metrics}
\author{Anton V. Solov'yov\thanks{Department of Theoretical Physics, Faculty of Physics, M. V. Lomonosov Moscow State University, Leninskie Gory, Moscow, 119991, Russia.\\
E-mail: a.v.solovyov@gmail.com}}
\date{}
\allowdisplaybreaks
\begin{document}
\maketitle
\begin{abstract}
We define natural volume forms on $n$-dimensional oriented pseudo-Finslerian manifolds with non-degenerate $m$-th root metrics. Our definitions of the natural volume forms depend on the parity of the positive integer $m>1$. The advantage of the stated definitions is their algebraic structure. The natural volume forms are expressed in terms of Cayley hyperdeterminants. In particular, the computation of the natural volume form does not require the difficult integration over the domain within the indicatrix in the tangent space $T_x M^n$ of the pseudo-Finslerian manifold at a point $x$.\\[10mm]
\end{abstract}
\section{Introduction}
The problem of defining a volume form on manifolds with additional geometric structures is not so trivial as it could seem at the first look. Although there is a general definition of volume forms on oriented differentiable manifolds, it contains some ambiguity. Indeed, let $x^1$,~\dots, $x^n\in\mathbb{R}$ be local coordinates in an $n$-dimensional oriented differentiable manifold $M^n$. By the definition~\cite{Kobayashi/Nomizu:book1996}, a differential $n$-form
\begin{equation}
\omega=\omega_{12\dots n}(x^1,\dots,x^n)dx^{1}\wedge\cdots\wedge dx^{n}
\label{eq:1}
\end{equation}
is called a \textit{volume form} (\textit{volume element}) if $\omega_{12\dots n}(x^1,\dots,x^n)>0$. Thus, \textit{any} positive function $\omega_{12\dots n}(x^1,\dots,x^n)$ generates the volume form~(\ref{eq:1}). However, there is no \textit{natural} choice of this function. Such a choice appears when $M^n$ has an additional geometric structure.

Let us recall  two classical examples. If $(M^{2n}, \Omega)$ is a $2n$-dimensional symplectic manifold with a closed non-degenerate differential 2-form $\Omega=\sum\limits_{a<b}\Omega_{ab}(x^1,\dots,x^{2n})dx^a\wedge dx^b$, then the natural volume form on it is $\omega=\underbrace{\Omega\wedge\cdots\wedge\Omega}_n$. If $(M^n, g)$ is an $n$-dimensional oriented Riemannian manifold with a metric tensor $g=g_{ab}(x^1,\dots,x^n)dx^a\otimes dx^b$, then the natural volume form on it is $\omega=\sqrt{\det[g_{ab}(x^1,\dots,x^n)]}dx^{1}\wedge\cdots\wedge dx^{n}$, where $g_{ab}(x^1,\dots,x^n)=g_{ba}(x^1,\dots,x^n)$. These examples illustrate a concept of a \textit{natural volume form}.

In Finslerian geometry~\cite{Rund:book1959,Bao/Chern/Shen:book2000,Chern/Shen:book2005}, we have several alternative definitions of the natural volume form. Let $(M^n, F)$ be an $n$-dimensional oriented Finslerian manifold with a line element $ds=F(x,dx)\geq 0$, where $F(x,\lambda dx)=\lambda F(x,dx)$ and $F(x,dx+dy)\leq F(x,dx)+F(x,dy)$ for any $x\in M^n$, $\lambda >0$, and $dx,dy\in T_x M^n$. If $\mathrm{Vol}(X)$ is the standard Euclidean volume of $X\subset\mathbb{R}^n$, then the natural volume form on $(M^n, F)$ is
\begin{equation}
\omega=\frac{\mathrm{Vol}(\text{the unit ball in $\mathbb{R}^n$})}{\mathrm{Vol}(\{v\in T_x M^n\mid F(x,v)\leq 1\})}dx^{1}\wedge\cdots\wedge dx^{n}.
\label{eq:2}
\end{equation}
This definition was proposed by H.~Busemann~\cite{Busemann:1947}. There are other definitions~\cite{Alvarez/Thompson:2004,Centore:1999} of the natural volume form on a Finslerian manifold. However, all of them essentially use properties of $T_x M^n$ as a normed space. In particular, $\mathrm{Vol}(\{v\in T_x M^n\mid F(x,v)\leq 1\})<\infty$.

In a \textit{pseudo-Finslerian manifold}, it is possible that $ds=F(x,dx)\ngeq 0$ and $F(x,dx+dy)\nleq F(x,dx)+F(x,dy)$ for some $dx,dy\in T_x M^n$. In other words, the tangent space $T_x M^n$ of a pseudo-Finslerian manifold at a point $x\in M^n$ is not a normed space. Therefore, $\mathrm{Vol}(\{v\in T_x M^n\mid F(x,v)\leq 1\})\nless\infty$ in some cases. The simplest example is the pseudo-Euclidean $xy$-plane, where $ds^2=dx^2 - dy^2$ and $\mathrm{Vol}(\{(x,y)\mid \sqrt{x^2 - y^2}\leq 1\})=\infty$. Thus, the definition~(\ref{eq:2}) \textit{fails} in pseudo-Finslerian manifolds.

In this paper, we define natural volume forms on $n$-dimensional oriented pseudo-Finslerian manifolds with the so-called ``$m$-th root metrics''
\begin{equation}
ds^m=g_{i_1 i_2\dots i_m}(x^1,\dots,x^n)dx^{i_1}dx^{i_2}\cdots dx^{i_m},
\label{eq:3}
\end{equation}
where $m>1$ is a positive integer and $i_1,i_2,\dots,i_m=1,2,\dots,n$. For $m=2$ and $m=4$, these metrics were considered by B.~Riemann in his famous lecture~\cite{Riemann:1868}. Differential geometry of Finslerian manifolds with the metrics~(\ref{eq:3}) was studied in~\cite{Matsumoto/Numata:1979,Shimada:1979}. There is a series of papers~\cite{Finkelstein:1986,Finkelstein/Holm:1986,Finkelstein/Holm:1987,Finkelstein:1988,Holm:1988a,Holm:1988b,Holm:1990,Honeycutt:1991,Solov'yov/Vladimirov:2001,Brody/Hughston:2005,Solov'yov:2006,Solov'yov:2011,Solov'yov:2015} on \textit{hyperspinors} (\textit{Finslerian N-spinors}) and their applications in which pseudo-Finslerian manifolds with the metrics~(\ref{eq:3}) appear intrinsically. It is evident that pseudo-Riemannian and pseudo-Finslerian manifolds are more adequate to relativistic physics than Riemannian and Finslerian ones. Classical and quantum field theories as well as relational approaches to quantum gravity require integration on manifolds. Therefore, any reasonable physical theory based on the metric~(\ref{eq:3}) will need a consistent definition of a volume form on the corresponding pseudo-Finslerian manifold.

An interesting definition of natural volume forms on pseudo-Finslerian manifolds with the metrics~(\ref{eq:3}) was proposed in~\cite{Kokarev:2004}. However, that definition is valid for even $m$ and special $n$ only. Below, we remove these restrictions by using a quite different construction of natural volume forms.
\section{General construction of the natural volume forms}
Let $(M^n, ds^m)$ be an $n$-dimensional oriented pseudo-Finslerian manifold with the metric~(\ref{eq:3}). For any volume form~(\ref{eq:1}), the function $\omega_{12\dots n}(x^1,\dots,x^n)$ is a component of a completely antisymmetric covariant tensor field in coordinates $x^1$,~\dots, $x^n$. Therefore, we have
\begin{equation}
\omega^\prime_{12\dots n}(x^{\prime 1},\dots,x^{\prime n})=\frac{\partial x^{j_1}}{\partial x^{\prime 1}}\frac{\partial x^{j_2}}{\partial x^{\prime 2}}\cdots\frac{\partial x^{j_n}}{\partial x^{\prime n}}\omega_{j_1 j_2\dots j_n}(x^1,\dots,x^n)
\label{eq:4}
\end{equation}
under coordinate transformations $x^j=x^j(x^{\prime 1},\dots,x^{\prime n})$, where all the indices run from 1 to $n$. Due to antisymmetry of the tensor components,
\begin{equation}
\omega_{j_1 j_2\dots j_n}(x^1,\dots,x^n)=\varepsilon_{j_1 j_2\dots j_n}\omega_{12\dots n}(x^1,\dots,x^n),
\label{eq:5}
\end{equation}
where $\varepsilon_{j_1 j_2\dots j_n}\equiv\varepsilon^{j_1 j_2\dots j_n}$ is the Levi-Civita symbol with $\varepsilon^{12\dots n}=1$. Inserting~(\ref{eq:5}) into~(\ref{eq:4}) and using the definition $\det[\frac{\partial x^j}{\partial x^{\prime i}}]=\varepsilon_{j_1 j_2\dots j_n}\frac{\partial x^{j_1}}{\partial x^{\prime 1}}\frac{\partial x^{j_2}}{\partial x^{\prime 2}}\cdots\frac{\partial x^{j_n}}{\partial x^{\prime n}}$ of the Jacobian determinant, we obtain the transformation law
\begin{equation}
\omega^\prime_{12\dots n}(x^{\prime 1},\dots,x^{\prime n})=\omega_{12\dots n}(x^1,\dots,x^n)\det\left[\frac{\partial x^j}{\partial x^{\prime i}}\right]
\label{eq:6}
\end{equation}
of the component $\omega_{12\dots n}>0$, where $\det[\frac{\partial x^j}{\partial x^{\prime i}}]>0$ because the manifold $M^n$ is oriented. Our aim is to express the function $\omega_{12\dots n}(x^1,\dots,x^n)$ in terms of the functions $g_{i_1 i_2\dots i_m}(x^1,\dots,x^n)$ from~(\ref{eq:3}) so that the transformation law~(\ref{eq:6}) is fulfilled for any $x\in M^n$.

The functions $g_{i_1 i_2\dots i_m}(x^1,\dots,x^n)$ are components of a completely symmetric pseudo-Finslerian metric tensor in  coordinates $x^1$,~\dots, $x^n$. Therefore,
\begin{equation}
g^\prime_{i_1 i_2\dots i_m}(x^{\prime 1},\dots,x^{\prime n})=\frac{\partial x^{j_1}}{\partial x^{\prime i_1}}\frac{\partial x^{j_2}}{\partial x^{\prime i_2}}\cdots\frac{\partial x^{j_m}}{\partial x^{\prime i_m}}g_{j_1 j_2\dots j_m}(x^1,\dots,x^n)
\label{eq:7}
\end{equation}
under coordinate transformations $x^j=x^j(x^{\prime 1},\dots,x^{\prime n})$. In order to deduce~(\ref{eq:6}) from~(\ref{eq:7}), we consider the following contraction
\begin{equation}
G^\prime_{i^1_1 i^1_2\dots i^1_n}=\varepsilon^{i^2_1 i^2_2\dots i^2_n}\cdots\varepsilon^{i^m_1 i^m_2\dots i^m_n}g^\prime_{i^1_1 i^2_1\dots i^m_1}g^\prime_{i^1_2 i^2_2\dots i^m_2}\cdots g^\prime_{i^1_n i^2_n\dots i^m_n}
\label{eq:8}
\end{equation}
(two-level indices are very useful in our computations and run from 1 to $n$ as well). Notice that $G^\prime_{i^1_1 i^1_2\dots i^1_n}$ are completely symmetric in their indices for odd $m>1$ and completely antisymmetric in their indices for even $m>0$. This fact is verified directly by using antisymmetry of the Levi-Civita symbol.

Let us insert~(\ref{eq:7}) into~(\ref{eq:8}). We obtain
\begin{multline}
G^\prime_{i^1_1\dots i^1_n}=\varepsilon^{i^2_1\dots i^2_n}\cdots\varepsilon^{i^m_1\dots i^m_n}g^\prime_{i^1_1 i^2_1\dots i^m_1}\cdots g^\prime_{i^1_n i^2_n\dots i^m_n}=\varepsilon^{i^2_1\dots i^2_n}\cdots\varepsilon^{i^m_1\dots i^m_n}\\
\times\left(\frac{\partial x^{j^1_1}}{\partial x^{\prime i^1_1}}\frac{\partial x^{j^2_1}}{\partial x^{\prime i^2_1}}\cdots\frac{\partial x^{j^m_1}}{\partial x^{\prime i^m_1}}g_{j^1_1 j^2_1\dots j^m_1}\right)
\cdots\left(\frac{\partial x^{j^1_n}}{\partial x^{\prime i^1_n}}\frac{\partial x^{j^2_n}}{\partial x^{\prime i^2_n}}\cdots\frac{\partial x^{j^m_n}}{\partial x^{\prime i^m_n}}g_{j^1_n j^2_n\dots j^m_n}\right)\\
=\frac{\partial x^{j^1_1}}{\partial x^{\prime i^1_1}}\cdots\frac{\partial x^{j^1_n}}{\partial x^{\prime i^1_n}}
\left(\varepsilon^{i^2_1\dots i^2_n}\frac{\partial x^{j^2_1}}{\partial x^{\prime i^2_1}}\cdots\frac{\partial x^{j^2_n}}{\partial x^{\prime i^2_n}}\right)\cdots\left(\varepsilon^{i^m_1\dots i^m_n}\frac{\partial x^{j^m_1}}{\partial x^{\prime i^m_1}}\cdots\frac{\partial x^{j^m_n}}{\partial x^{\prime i^m_n}}\right)\\
\times g_{j^1_1 j^2_1\dots j^m_1}\cdots g_{j^1_n j^2_n\dots j^m_n}=\frac{\partial x^{j^1_1}}{\partial x^{\prime i^1_1}}\cdots\frac{\partial x^{j^1_n}}{\partial x^{\prime i^1_n}}\varepsilon^{j^2_1\dots j^2_n}\det\left[\frac{\partial x^j}{\partial x^{\prime i}}\right]\cdots\varepsilon^{j^m_1\dots j^m_n}\\
\times\det\left[\frac{\partial x^j}{\partial x^{\prime i}}\right]g_{j^1_1 j^2_1\dots j^m_1}\cdots g_{j^1_n j^2_n\dots j^m_n}=\left(\det\left[\frac{\partial x^j}{\partial x^{\prime i}}\right]\right)^{m-1}\frac{\partial x^{j^1_1}}{\partial x^{\prime i^1_1}}\cdots\frac{\partial x^{j^1_n}}{\partial x^{\prime i^1_n}}\\
\times\varepsilon^{j^2_1\dots j^2_n}\cdots\varepsilon^{j^m_1\dots j^m_n}g_{j^1_1 j^2_1\dots j^m_1}\cdots g_{j^1_n j^2_n\dots j^m_n}.
\label{eq:9}
\end{multline}
Using the notation
\begin{equation}
G_{j^1_1\dots j^1_n}=\varepsilon^{j^2_1\dots j^2_n}\cdots\varepsilon^{j^m_1\dots j^m_n}g_{j^1_1 j^2_1\dots j^m_1}\cdots g_{j^1_n j^2_n\dots j^m_n},
\label{eq:10}
\end{equation}
we can rewrite~(\ref{eq:9}) as
\begin{equation}
G^\prime_{i^1_1\dots i^1_n}=\left(\det\left[\frac{\partial x^j}{\partial x^{\prime i}}\right]\right)^{m-1}\frac{\partial x^{j^1_1}}{\partial x^{\prime i^1_1}}\cdots\frac{\partial x^{j^1_n}}{\partial x^{\prime i^1_n}}G_{j^1_1\dots j^1_n}.
\label{eq:11}
\end{equation}
Notice the identical structure of the expressions~(\ref{eq:8}) and~(\ref{eq:10}). Thus,~(\ref{eq:11}) is the transformation law of a tensor density of weight $m-1$.
\subsection{The natural volume forms for even $m>0$}
Let us consider another contraction
\begin{equation}
G^\prime=\varepsilon^{i^1_1 i^1_2\dots i^1_n}G^\prime_{i^1_1 i^1_2\dots i^1_n}.
\label{eq:12}
\end{equation}
It is evident that $G^\prime\equiv 0$ for odd $m>1$ ($\varepsilon^{i^1_1 i^1_2\dots i^1_n}$ are antisymmetric and $G^\prime_{i^1_1 i^1_2\dots i^1_n}$ are symmetric in this case). Therefore, we focus on even $m>0$.

Inserting~(\ref{eq:11}) into~(\ref{eq:12}), we have
\begin{multline}
G^\prime=\varepsilon^{i^1_1\dots i^1_n}G^\prime_{i^1_1\dots i^1_n}=\varepsilon^{i^1_1\dots i^1_n}\left(\det\left[\frac{\partial x^j}{\partial x^{\prime i}}\right]\right)^{m-1}\frac{\partial x^{j^1_1}}{\partial x^{\prime i^1_1}}\cdots\frac{\partial x^{j^1_n}}{\partial x^{\prime i^1_n}}G_{j^1_1\dots j^1_n}\\
=\left(\det\left[\frac{\partial x^j}{\partial x^{\prime i}}\right]\right)^{m-1}\varepsilon^{i^1_1\dots i^1_n}\frac{\partial x^{j^1_1}}{\partial x^{\prime i^1_1}}\cdots\frac{\partial x^{j^1_n}}{\partial x^{\prime i^1_n}}G_{j^1_1\dots j^1_n}=\left(\det\left[\frac{\partial x^j}{\partial x^{\prime i}}\right]\right)^{m-1}\\
\times\varepsilon^{j^1_1\dots j^1_n}\det\left[\frac{\partial x^j}{\partial x^{\prime i}}\right]G_{j^1_1\dots j^1_n}=\left(\det\left[\frac{\partial x^j}{\partial x^{\prime i}}\right]\right)^m \varepsilon^{j^1_1\dots j^1_n}G_{j^1_1\dots j^1_n}.
\label{eq:13}
\end{multline}
Using the notation
\begin{equation}
G=\varepsilon^{j^1_1\dots j^1_n}G_{j^1_1\dots j^1_n},
\label{eq:14}
\end{equation}
we can rewrite~(\ref{eq:13}) as
\begin{equation}
G^\prime=\left(\det\left[\frac{\partial x^j}{\partial x^{\prime i}}\right]\right)^m G.
\label{eq:15}
\end{equation}
It is easy to see that~(\ref{eq:15}) is  the transformation law of a scalar density of weight $m$.

Let us compare~(\ref{eq:6}) with~(\ref{eq:15}). It is clear that $\sqrt[m]{|G|}$ is a suitable candidate for the function $\omega_{12\dots n}(x^1,\dots,x^n)$. However, $G^\prime_{i^1_1 i^1_2\dots i^1_n}$ and $G_{j^1_1 j^1_2\dots j^1_n}$ are completely antisymmetric for even $m>0$. Therefore,
\begin{equation}
G^\prime_{i^1_1 i^1_2\dots i^1_n}=\varepsilon_{i^1_1 i^1_2\dots i^1_n}G^\prime_{1 2\dots n}
\quad\text{and}\quad
G_{j^1_1 j^1_2\dots j^1_n}=\varepsilon_{j^1_1 j^1_2\dots j^1_n}G_{1 2\dots n}.
\label{eq:16}
\end{equation}
Inserting~(\ref{eq:16}) into~(\ref{eq:12}) and~(\ref{eq:14}), we obtain
\begin{equation}
G^\prime=n!\,G^\prime_{1 2\dots n}\quad\text{and}\quad G=n!\,G_{1 2\dots n},
\label{eq:17}
\end{equation}
where the known formula $\varepsilon^{i^1_1 i^1_2\dots i^1_n}\varepsilon_{i^1_1 i^1_2\dots i^1_n}=\varepsilon^{j^1_1 j^1_2\dots j^1_n}\varepsilon_{j^1_1 j^1_2\dots j^1_n}=n!$ has been used. Thus,~(\ref{eq:15}) and~(\ref{eq:17}) imply
\begin{equation}
|G^\prime_{1 2\dots n}|^{1/m}=|G_{1 2\dots n}|^{1/m}\det\left[\frac{\partial x^j}{\partial x^{\prime i}}\right].
\label{eq:18}
\end{equation}
If $\omega_{12\dots n}(x^1,\dots,x^n)=|G_{1 2\dots n}|^{1/m}$, then~(\ref{eq:6}) and~(\ref{eq:18}) coincide.

We will use the following notation
\begin{multline}
\mathrm{hdet}[g_{i_1 i_2\dots i_m}(x^1,\dots,x^n)]\equiv G_{1 2\dots n}\\
=\varepsilon^{j^2_1 j^2_2\dots j^2_n}\cdots\varepsilon^{j^m_1 j^m_2\dots j^m_n}g_{1 j^2_1\dots j^m_1}g_{2 j^2_2\dots j^m_2}\cdots g_{n j^2_n\dots j^m_n}.
\label{eq:19}
\end{multline}
It is interesting that homogeneous polynomials of type~(\ref{eq:19}) were introduced by A.~Cayley (in other terms and with respect to independent variables, not functions) in the second part of the paper~\cite{Cayley:1843}. Later on, he called them \textit{hyperdeterminants}~\cite{Cayley:1845}. An excellent exposition of hyperdeterminants and their modifications can be found in the book~\cite{Sokolov:book1960}. It is evident that the hyperdeterminant~(\ref{eq:19}) becomes the usual determinant for $m=2$: $\mathrm{hdet}[g_{i_1 i_2}]=\det[g_{i_1 i_2}]$.

Using the notation~(\ref{eq:19}), we rewrite~(\ref{eq:18}) finally as
\begin{equation}
\bigl|\mathrm{hdet}[g^\prime_{i_1 i_2\dots i_m}]\bigr|^{1/m}=\bigl|\mathrm{hdet}[g_{i_1 i_2\dots i_m}]\bigr|^{1/m}\det\left[\frac{\partial x^j}{\partial x^{\prime i}}\right].
\label{eq:20}
\end{equation}
Because of~(\ref{eq:1}), (\ref{eq:6}), and~(\ref{eq:20}), we propose
\begin{equation}
\omega=\bigl|\mathrm{hdet}[g_{i_1 i_2\dots i_m}(x^1,\dots,x^n)]\bigr|^{1/m}dx^{1}\wedge\cdots\wedge dx^{n}
\label{eq:21}
\end{equation}
as the natural volume form on the $n$-dimensional oriented pseudo-Finslerian manifold $(M^n, ds^m)$ with the metric~(\ref{eq:3}) for even $m>0$. Of course, the hyperdeterminant~(\ref{eq:19}) can vanish accidentally. In order to avoid such situations, we call the metric~(\ref{eq:3}) \textit{non-degenerate} if
\begin{equation}
\mathrm{hdet}[g_{i_1 i_2\dots i_m}(x^1,\dots,x^n)]\neq 0
\label{eq:22}
\end{equation}
for any $x\in M^n$ and even $m>0$. Thus,~(\ref{eq:21}) is the natural volume form on $(M^n, ds^m)$ with the non-degenerate metric~(\ref{eq:3}) for which the requirement~(\ref{eq:22}) is obligatory and $m>0$ is even. 

In the case of $m=2$,~(\ref{eq:21}) is the well-known pseudo-Riemannian natural volume form $\omega=\bigl|\det[g_{i_1 i_2}(x^1,\dots,x^n)]\bigr|^{1/2}dx^{1}\wedge\cdots\wedge dx^{n}$. Let us consider the more interesting example when $m=4$ and $n=2$. In this case, the metric~(\ref{eq:3}) is
\begin{equation}
ds^4=g_{i_1 i_2 i_3 i_4}(x^1,x^2)dx^{i_1}dx^{i_2}dx^{i_3} dx^{i_4},
\label{eq:23}
\end{equation}
where $i_1$, $i_2$, $i_3$, $i_4=1,2$. By computing the hyperdeterminant~(\ref{eq:19}) for the metric~(\ref{eq:23}), we obtain
\begin{multline}
\mathrm{hdet}[g_{i_1 i_2 i_3 i_4}(x^1,x^2)]=\varepsilon^{j^2_1 j^2_2}\varepsilon^{j^3_1 j^3_2}\varepsilon^{j^4_1 j^4_2}g_{1 j^2_1 j^3_1 j^4_1}g_{2 j^2_2 j^3_2 j^4_2}\\
=g_{1111}g_{2222}-g_{1112}g_{2221}-g_{1121}g_{2212}+g_{1122}g_{2211}\\
-g_{1211}g_{2122}+g_{1212}g_{2121}+g_{1221}g_{2112}-g_{1222}g_{2111}.
\label{eq:24}
\end{multline}
However, $g_{i_1 i_2 i_3 i_4}$ are symmetric in all the indices. Therefore,
\begin{align}
g_{1112}&=g_{1121}=g_{1211}=g_{2111},\nonumber\\
g_{1122}&=g_{1221}=g_{2121}=g_{1212}=g_{2112}=g_{2211},\nonumber\\
g_{1222}&=g_{2221}=g_{2212}=g_{2122}.
\label{eq:25}
\end{align}
Because of~(\ref{eq:25}), we can choose $g_{1111}$, $g_{1112}$, $g_{1122}$, $g_{1222}$, and $g_{2222}$ as independent components of the pseudo-Finslerian metric tensor. In this case,~(\ref{eq:24}) has the form
\begin{equation}
\mathrm{hdet}[g_{i_1 i_2 i_3 i_4}(x^1,x^2)]=g_{1111}g_{2222}-4g_{1112}g_{1222}+3(g_{1122})^2.
\label{eq:26}
\end{equation}
Inserting~(\ref{eq:26}) into~(\ref{eq:21}) and~(\ref{eq:22}), we see that the natural volume form on $(M^2, ds^4)$ is
\begin{equation}
\omega=\bigl|g_{1111}g_{2222}-4g_{1112}g_{1222}+3(g_{1122})^2\bigr|^{1/4}dx^{1}\wedge dx^{2}
\label{eq:27}
\end{equation}
and the metric~(\ref{eq:23}) is non-degenerate if
\begin{equation}
g_{1111}g_{2222}-4g_{1112}g_{1222}+3(g_{1122})^2\neq 0
\label{eq:28}
\end{equation}
for any $x\in M^2$. The formulas~(\ref{eq:27}) and~(\ref{eq:28}) give us a non-trivial example of the pseudo-Finslerian natural volume form for even $m>0$.
\subsection{The natural volume forms for odd $m>1$}
Unfortunately, the definition~(\ref{eq:21}) is not suitable for odd $m>1$. Indeed, we mentioned above that~(\ref{eq:12}) and, of course,~(\ref{eq:14}) vanish identically in this case. Moreover, the transformation law~(\ref{eq:20}) fails for odd $m>1$ so that the hyperdeterminant~(\ref{eq:19}) cannot be used directly to define the natural volume form on $(M^n, ds^m)$. However, we can iterate the above construction in the following way.

Notice that the transformation laws~(\ref{eq:7}) and~(\ref{eq:11}) differ only in the multiplier $\left(\det\left[\frac{\partial x^j}{\partial x^{\prime i}}\right]\right)^{m-1}$. Let us replace $g^\prime_{i^1_1 i^2_1\dots i^m_1}$ with $G^\prime_{i^1_1 i^2_1\dots i^n_1}$ everywhere in~(\ref{eq:8}), i.e., consider the contraction
\begin{equation}
\widetilde{G}^\prime_{i^1_1 i^1_2\dots i^1_n}=\varepsilon^{i^2_1 i^2_2\dots i^2_n}\cdots\varepsilon^{i^n_1 i^n_2\dots i^n_n}G^\prime_{i^1_1 i^2_1\dots i^n_1}G^\prime_{i^1_2 i^2_2\dots i^n_2}\cdots G^\prime_{i^1_n i^2_n\dots i^n_n}.
\label{eq:29}
\end{equation}
Inserting~(\ref{eq:11}) into~(\ref{eq:29}), we obtain
\begin{multline}
\widetilde{G}^\prime_{i^1_1\dots i^1_n}=\varepsilon^{i^2_1\dots i^2_n}\cdots\varepsilon^{i^n_1\dots i^n_n}G^\prime_{i^1_1 i^2_1\dots i^n_1}\cdots G^\prime_{i^1_n i^2_n\dots i^n_n}=\varepsilon^{i^2_1\dots i^2_n}\cdots\varepsilon^{i^n_1\dots i^n_n}\\
\times\left(\det\left[\frac{\partial x^j}{\partial x^{\prime i}}\right]\right)^{m-1}
\left(\frac{\partial x^{j^1_1}}{\partial x^{\prime i^1_1}}\frac{\partial x^{j^2_1}}{\partial x^{\prime i^2_1}}\cdots\frac{\partial x^{j^n_1}}{\partial x^{\prime i^n_1}}G_{j^1_1 j^2_1\dots j^n_1}\right)\cdots\\
\times\left(\det\left[\frac{\partial x^j}{\partial x^{\prime i}}\right]\right)^{m-1}\left(\frac{\partial x^{j^1_n}}{\partial x^{\prime i^1_n}}\frac{\partial x^{j^2_n}}{\partial x^{\prime i^2_n}}\cdots\frac{\partial x^{j^n_n}}{\partial x^{\prime i^n_n}}G_{j^1_n j^2_n\dots j^n_n}\right)\\
=\left(\det\left[\frac{\partial x^j}{\partial x^{\prime i}}\right]\right)^{(m-1)n}\frac{\partial x^{j^1_1}}{\partial x^{\prime i^1_1}}\cdots\frac{\partial x^{j^1_n}}{\partial x^{\prime i^1_n}}\left(\varepsilon^{i^2_1\dots i^2_n}\frac{\partial x^{j^2_1}}{\partial x^{\prime i^2_1}}\cdots\frac{\partial x^{j^2_n}}{\partial x^{\prime i^2_n}}\right)\cdots\\
\times\left(\varepsilon^{i^n_1\dots i^n_n}\frac{\partial x^{j^n_1}}{\partial x^{\prime i^n_1}}\cdots\frac{\partial x^{j^n_n}}{\partial x^{\prime i^n_n}}\right)G_{j^1_1 j^2_1\dots j^n_1}\cdots G_{j^1_n j^2_n\dots j^n_n}=\left(\det\left[\frac{\partial x^j}{\partial x^{\prime i}}\right]\right)^{(m-1)n}\\
\times\frac{\partial x^{j^1_1}}{\partial x^{\prime i^1_1}}\cdots\frac{\partial x^{j^1_n}}{\partial x^{\prime i^1_n}}\varepsilon^{j^2_1\dots j^2_n}\det\left[\frac{\partial x^j}{\partial x^{\prime i}}\right]\cdots\varepsilon^{j^n_1\dots j^n_n}
\det\left[\frac{\partial x^j}{\partial x^{\prime i}}\right]\\
\times G_{j^1_1 j^2_1\dots j^n_1}\cdots G_{j^1_n j^2_n\dots j^n_n}=\left(\det\left[\frac{\partial x^j}{\partial x^{\prime i}}\right]\right)^{mn-1}\frac{\partial x^{j^1_1}}{\partial x^{\prime i^1_1}}\cdots\frac{\partial x^{j^1_n}}{\partial x^{\prime i^1_n}}\\
\times\varepsilon^{j^2_1\dots j^2_n}\cdots\varepsilon^{j^n_1\dots j^n_n}G_{j^1_1 j^2_1\dots j^n_1}\cdots G_{j^1_n j^2_n\dots j^n_n}.
\label{eq:30}
\end{multline}
By using the notation
\begin{equation}
\widetilde{G}_{j^1_1\dots j^1_n}=\varepsilon^{j^2_1\dots j^2_n}\cdots\varepsilon^{j^n_1\dots j^n_n}G_{j^1_1 j^2_1\dots j^n_1}\cdots G_{j^1_n j^2_n\dots j^n_n},
\label{eq:31}
\end{equation}
we rewrite~(\ref{eq:30}) in the following form
\begin{equation}
\widetilde{G}^\prime_{i^1_1\dots i^1_n}=\left(\det\left[\frac{\partial x^j}{\partial x^{\prime i}}\right]\right)^{mn-1}\frac{\partial x^{j^1_1}}{\partial x^{\prime i^1_1}}\cdots\frac{\partial x^{j^1_n}}{\partial x^{\prime i^1_n}}\widetilde{G}_{j^1_1\dots j^1_n}.
\label{eq:32}
\end{equation}
Thus,~(\ref{eq:32}) is the transformation law of a tensor density of weight $mn-1$.

In order to construct a scalar density from $\widetilde{G}^\prime_{i^1_1\dots i^1_n}$, we compute the total contraction
\begin{equation}
\widetilde{G}^\prime=\varepsilon^{i^1_1 i^1_2\dots i^1_n}\widetilde{G}^\prime_{i^1_1 i^1_2\dots i^1_n}.
\label{eq:33}
\end{equation}
It is evident that $\widetilde{G}^\prime\equiv 0$ for odd $n>1$ ($\varepsilon^{i^1_1 i^1_2\dots i^1_n}$ are antisymmetric and $\widetilde{G}^\prime_{i^1_1 i^1_2\dots i^1_n}$ are symmetric in this case). Therefore, we are forced to consider \textit{even $n>0$ only}.

Inserting~(\ref{eq:32}) into~(\ref{eq:33}), we have
\begin{multline}
\widetilde{G}^\prime=\varepsilon^{i^1_1\dots i^1_n}\widetilde{G}^\prime_{i^1_1\dots i^1_n}=\varepsilon^{i^1_1\dots i^1_n}\left(\det\left[\frac{\partial x^j}{\partial x^{\prime i}}\right]\right)^{mn-1}\frac{\partial x^{j^1_1}}{\partial x^{\prime i^1_1}}\cdots\frac{\partial x^{j^1_n}}{\partial x^{\prime i^1_n}}\widetilde{G}_{j^1_1\dots j^1_n}\\
=\left(\det\left[\frac{\partial x^j}{\partial x^{\prime i}}\right]\right)^{mn-1}\varepsilon^{i^1_1\dots i^1_n}\frac{\partial x^{j^1_1}}{\partial x^{\prime i^1_1}}\cdots\frac{\partial x^{j^1_n}}{\partial x^{\prime i^1_n}}\widetilde{G}_{j^1_1\dots j^1_n}=\left(\det\left[\frac{\partial x^j}{\partial x^{\prime i}}\right]\right)^{mn-1}\\
\times\varepsilon^{j^1_1\dots j^1_n}\det\left[\frac{\partial x^j}{\partial x^{\prime i}}\right]\widetilde{G}_{j^1_1\dots j^1_n}=\left(\det\left[\frac{\partial x^j}{\partial x^{\prime i}}\right]\right)^{mn}\varepsilon^{j^1_1\dots j^1_n}\widetilde{G}_{j^1_1\dots j^1_n}.
\label{eq:34}
\end{multline}
By using the notation
\begin{equation}
\widetilde{G}=\varepsilon^{j^1_1\dots j^1_n}\widetilde{G}_{j^1_1\dots j^1_n},
\label{eq:35}
\end{equation}
we can rewrite~(\ref{eq:34}) as
\begin{equation}
\widetilde{G}^\prime=\left(\det\left[\frac{\partial x^j}{\partial x^{\prime i}}\right]\right)^{mn}\widetilde{G}.
\label{eq:36}
\end{equation}
It is easy to see that~(\ref{eq:36}) is  the transformation law of a scalar density of weight $mn$.

Because of~(\ref{eq:29}) and~(\ref{eq:31}), the components $\widetilde{G}^\prime_{i^1_1 i^1_2\dots i^1_n}$ and $\widetilde{G}_{j^1_1 j^1_2\dots j^1_n}$ are completely antisymmetric for even $n>0$. Therefore,
\begin{equation}
\widetilde{G}^\prime_{i^1_1 i^1_2\dots i^1_n}=\varepsilon_{i^1_1 i^1_2\dots i^1_n}\widetilde{G}^\prime_{1 2\dots n}
\quad\text{and}\quad
\widetilde{G}_{j^1_1 j^1_2\dots j^1_n}=\varepsilon_{j^1_1 j^1_2\dots j^1_n}\widetilde{G}_{1 2\dots n}.
\label{eq:37}
\end{equation}
By inserting~(\ref{eq:37}) into~(\ref{eq:33}) and~(\ref{eq:35}), we obtain
\begin{equation}
\widetilde{G}^\prime=n!\,\widetilde{G}^\prime_{1 2\dots n}\quad\text{and}\quad \widetilde{G}=n!\,\widetilde{G}_{1 2\dots n}.
\label{eq:38}
\end{equation}
At the same time, $\widetilde{G}_{1 2\dots n}$ is the hyperdeterminant
\begin{multline}
\mathrm{hdet}[G_{i_1 i_2\dots i_n}(x^1,\dots,x^n)]\equiv\widetilde{G}_{1 2\dots n}\\
=\varepsilon^{j^2_1 j^2_2\dots j^2_n}\cdots\varepsilon^{j^n_1 j^n_2\dots j^n_n}G_{1 j^2_1\dots j^n_1}G_{2 j^2_2\dots j^n_2}\cdots G_{n j^2_n\dots j^n_n}
\label{eq:39}
\end{multline}
(compare~(\ref{eq:39}) with~(\ref{eq:19}) for clarity). Thus,~(\ref{eq:36}), (\ref{eq:38}), and~(\ref{eq:39}) imply
\begin{equation}
\bigl|\mathrm{hdet}[G^\prime_{i_1 i_2\dots i_n}]\bigr|^{1/(mn)}=\bigl|\mathrm{hdet}[G_{i_1 i_2\dots i_n}]\bigr|^{1/(mn)}\det\left[\frac{\partial x^j}{\partial x^{\prime i}}\right].
\label{eq:40}
\end{equation}

Because of~(\ref{eq:1}), (\ref{eq:6}), and~(\ref{eq:40}), we propose
\begin{equation}
\omega=\bigl|\mathrm{hdet}[G_{i_1 i_2\dots i_n}(x^1,\dots,x^n)]\bigr|^{1/(mn)}dx^{1}\wedge\cdots\wedge dx^{n}
\label{eq:41}
\end{equation}
as the natural volume form on the $n$-dimensional oriented pseudo-Finslerian manifold $(M^n, ds^m)$ with the metric~(\ref{eq:3}) for odd  $m>1$ and even $n>0$ (notice that the functions $G_{j^1_1 j^1_2\dots j^1_n}$ defined by~(\ref{eq:10}) depend on~$m$ too). Again, the hyperdeterminant~(\ref{eq:39}) can vanish in some cases. In order to avoid such situations, we call the metric~(\ref{eq:3}) \textit{non-degenerate} if
\begin{equation}
\mathrm{hdet}[G_{i_1 i_2\dots i_n}(x^1,\dots,x^n)]\neq 0
\label{eq:42}
\end{equation}
for any $x\in M^n$, odd  $m>1$, and even $n>0$. Thus,~(\ref{eq:41}) is the natural volume form on $(M^n, ds^m)$ with the non-degenerate metric~(\ref{eq:3}) for which the requirement~(\ref{eq:42}) is obligatory, $m>1$ is odd, and $n>0$ is even.

Let us consider the simplest example when $m=3$ and $n=2$. In this case, the metric~(\ref{eq:3}) is
\begin{equation}
ds^3=g_{i_1 i_2 i_3}(x^1,x^2)dx^{i_1}dx^{i_2}dx^{i_3},
\label{eq:43}
\end{equation}
where $i_1$, $i_2$, $i_3=1,2$. By computing the hyperdeterminant~(\ref{eq:39}) with the help of the definition~(\ref{eq:10}) for the metric~(\ref{eq:43}), we obtain
\begin{multline}
\mathrm{hdet}[G_{i_1 i_2}(x^1,x^2)]=\det[G_{i_1 i_2}(x^1,x^2)]=\varepsilon^{j^2_1 j^2_2}G_{1 j^2_1}G_{2 j^2_2}\\
=\varepsilon^{j^2_1 j^2_2}\varepsilon^{k^2_1 k^2_2}\varepsilon^{k^3_1 k^3_2}g_{1 k^2_1 k^3_1}g_{j^2_1 k^2_2 k^3_2}\varepsilon^{l^2_1 l^2_2}\varepsilon^{l^3_1 l^3_2}g_{2 l^2_1 l^3_1}g_{j^2_2 l^2_2 l^3_2}\\
=-(g_{111}g_{222})^2 -(g_{112}g_{221})^2 -(g_{121}g_{212})^2 -(g_{122}g_{211})^2\\
+2(g_{111}g_{112}g_{221}g_{222}+g_{111}g_{121}g_{212}g_{222}+g_{111}g_{122}g_{211}g_{222}\\
+g_{112}g_{121}g_{212}g_{221}+g_{112}g_{122}g_{211}g_{221}+g_{121}g_{122}g_{211}g_{212})\\
-4(g_{111}g_{122}g_{212}g_{221}+g_{112}g_{121}g_{211}g_{222}).
\label{eq:44}
\end{multline}
It is surprising that~(\ref{eq:44}) multiplied by $-1$ coincides with one of the simplest \textit{complete hyperdeterminants} introduced by A.~Cayley (with respect to independent variables, not functions) in the paper~\cite{Cayley:1845}. This  hyperdeterminant and its generalizations are discriminants of certain multilinear forms. The discriminant aspects of hyperdeterminants are studied in the fundamental book~\cite{Gelfand/Kapranov/Zelevinsky:book1994}.

Let us take into account that $g_{i_1 i_2 i_3}$ are symmetric in all the indices. Therefore,
\begin{equation}
g_{112}=g_{121}=g_{211},\quad g_{122}=g_{212}=g_{221}.
\label{eq:45}
\end{equation}
Because of~(\ref{eq:45}), we can choose $g_{111}$, $g_{112}$, $g_{122}$, and $g_{222}$ as independent components of the pseudo-Finslerian metric tensor. In this case,~(\ref{eq:44}) has the form
\begin{multline}
\mathrm{hdet}[G_{i_1 i_2}(x^1,x^2)]=-(g_{111}g_{222})^2 + 6g_{111}g_{112}g_{122}g_{222}\\
- 4g_{111}(g_{122})^3 - 4(g_{112})^3g_{222}+3(g_{112}g_{122})^2.
\label{eq:46}
\end{multline}
Inserting~(\ref{eq:46}) into~(\ref{eq:41}) and~(\ref{eq:42}), we see that the natural volume form on $(M^2, ds^3)$ is
\begin{multline}
\omega=\bigl|-(g_{111}g_{222})^2 + 6g_{111}g_{112}g_{122}g_{222} - 4g_{111}(g_{122})^3\\
- 4(g_{112})^3g_{222}+3(g_{112}g_{122})^2 \bigr|^{1/6}dx^{1}\wedge dx^{2}
\label{eq:47}
\end{multline}
and the metric~(\ref{eq:43}) is non-degenerate if
\begin{multline}
-(g_{111}g_{222})^2 + 6g_{111}g_{112}g_{122}g_{222} - 4g_{111}(g_{122})^3\\
- 4(g_{112})^3g_{222}+3(g_{112}g_{122})^2\neq 0
\label{eq:48}
\end{multline}
for any $x\in M^2$. The formulas~(\ref{eq:47}) and~(\ref{eq:48}) give us a non-trivial example of the pseudo-Finslerian natural volume form for odd $m>1$.
\section{Conclusion}
In this paper, we have constructed the natural volume forms on the $n$-dimensional oriented pseudo-Finslerian manifolds $(M^n, ds^m)$ with the non-degenerate $m$-th root metrics~(\ref{eq:3}).

It is important that our definition of the natural volume form depends on the \textit{parity} of the positive integer $m$. For even $m>0$ and integer $n>1$, the natural volume form $\omega$ on $(M^n, ds^m)$ is defined by the formula~(\ref{eq:21}). If $m=2$, then~(\ref{eq:21}) becomes the standard pseudo-Riemannian natural volume form. For odd $m>1$ and even $n>0$, the natural volume form $\omega$ on $(M^n, ds^m)$ is defined by the formula~(\ref{eq:41}). The condition that the metric~(\ref{eq:3}) is non-degenerate depends also on the parity of $m$ and is given by~(\ref{eq:22}) or~(\ref{eq:42}).

Unfortunately, the author was unable to obtain the general expression of  the natural volume form $\omega$ on $(M^n, ds^m)$ for odd $m>1$ and odd $n>1$ (although there are several special examples of $\omega$ in this case). Moreover, the author does not assert that~(\ref{eq:21}) and~(\ref{eq:41}) are the unique definitions of the natural volume forms on $(M^n, ds^m)$, but he believes that~(\ref{eq:21}) and~(\ref{eq:41}) are the simplest ones.

The advantage of the definitions~(\ref{eq:21}) and~(\ref{eq:41}) is their algebraic structure. According to~(\ref{eq:21}) and~(\ref{eq:41}), the natural volume form $\omega$ on $(M^n, ds^m)$ is expressed in terms of the Cayley hyperdeterminants. In particular, the computation of $\omega$ does not require the difficult integration over the domain within the indicatrix of the tangent space $T_x M^n$ as in~(\ref{eq:2}).

The definitions~(\ref{eq:21}) and~(\ref{eq:41}), up to constant factors, were invented respectively in~\cite{Solov'yov:tez2012} and~\cite{Solov'yov:tez2013}, but at that time, the author did not understand deep interconnections between them and different Cayley hyperdeterminants.

\section*{Acknowledgements}

The author is grateful to organizers of the VIIIth International Conference ``Finsler Extensions of Relativity Theory'' for hospitality and to its participants for helpful discussions of different aspects of this paper. I especially thank R.~R.~Aidagulov, V.~Balan, G.~I.~Garas'ko, A.~V.~Koganov, S.~S.~Kokarev, D.~G.~Pavlov, and Shamil Shakirov. Finally, I thank my father, Vasiliy Ivanovich Solov'yov, for all\dots

\bibliographystyle{ieeetr}
\bibliography{volume_forms}
\end{document}